\documentclass[reqno, 12pt]{amsart}
 \usepackage[a4paper, hmargin = 1.25 in, vmargin = 1.5 in]{geometry}
 \usepackage{amsmath, amssymb, amsfonts}
 \usepackage{amsthm}
 \usepackage{xcolor}
 \usepackage{graphicx}

 \newtheorem{theorem}{Theorem}[section]
 
 \newtheorem{exe}{Example}
 \newtheorem{lemma}[theorem]{Lemma}
 
 \theoremstyle{definition}
 \newtheorem{definition}{Definition}[section]
 \theoremstyle{remark}
 
 \numberwithin{equation}{section}

\begin{document}

\begin{center}

\begin{title}
\title{\bf\Large{{Numerical Computation of Exponential Functions of Nabla Fractional Calculus}}}
\end{title}

\vskip 0.25 in

\begin{author}
\author {Jagan Mohan Jonnalagadda\footnote[1]{Department of Mathematics, Birla Institute of Technology and Science Pilani, Hyderabad - 500078, Telangana, India. email: {j.jaganmohan@hotmail.com}}}
\end{author}

\end{center}

\vskip 0.25 in

\noindent{\bf Abstract:} In this article, we illustrate the asymptotic behaviour of exponential functions of nabla fractional calculus. For this purpose, we propose a novel matrix technique to compute these functions numerically.

\vskip 0.25 in


\vskip 0.25 in


\vskip 0.25 in

\section{Introduction \& Preliminaries}
Nabla fractional calculus is an integrated theory of arbitrary order sums and differences. The concept of nabla fractional difference traces back to the works of Miller \& Ross \cite{Mi}, Gray \& Zhang \cite{Gray}, Atici \& Eloe \cite{At 1}, and Anastassiou \cite{An}. During the past one decade, there has been an increasing interest in this field. For a detailed introduction on the evolution of nabla fractional calculus, we refer to \cite{Go} and the references therein.

We use the following notations, definitions and known results of nabla fractional calculus throughout the article. Denote by $\mathbb{N}_{a} = \{a, a + 1, a + 2, \ldots\}$ and $\mathbb{N}^{b}_{a} = \{a, a + 1, a + 2, \ldots, b\}$ for any $a$, $b \in \mathbb{R}$ such that $b - a \in \mathbb{N}_{1}$. The backward jump operator $\rho : \mathbb{N}_{a + 1} \rightarrow \mathbb{N}_{a}$ is defined by $$\rho(t) = t - 1, \quad t \in \mathbb{N}_{a + 1}.$$ Define the $\mu^{th}$-order nabla fractional Taylor monomial by $$H_{\mu}(t, a) = \frac{\Gamma(t - a + \mu)}{\Gamma(t - a)\Gamma(\mu + 1)}, \quad \mu \in \mathbb{R} \setminus \{\ldots, -2, -1\},$$ provided the right-hand side of this equation is sensible. Here $\Gamma(\cdot)$ denotes the Euler gamma function.

\begin{lemma} \cite{Go}
We observe the following properties of nabla fractional Taylor monomials.
\begin{enumerate}
\item $H_{\mu}(t, a) = 0$ for all $\mu \in \{\ldots, -2, -1\}$ and $t \in \mathbb{N}_{a}$. 
\item $H_{\mu}(t, \rho(t)) = 1$ for all $\mu \in \mathbb{R} \setminus \{\ldots, -2, -1\}$ and $t \in \mathbb{N}_{a}$. 
\item $H_{\mu}(t, t) = 0$ for all $\mu \in \mathbb{R} \setminus \{\ldots, -2, -1\}$ and $t \in \mathbb{N}_{a}$. 
\end{enumerate}
\end{lemma}

\begin{definition} \cite{Bo}
Let $u: \mathbb{N}_{a} \rightarrow \mathbb{R}$. The first order backward (nabla) difference of $u$ is defined by $$\big{(}\nabla u\big{)}(t) = u(t) - u(t - 1), \quad t \in \mathbb{N}_{a + 1}.$$
\end{definition}

\begin{definition} \cite{Go}
Let $u: \mathbb{N}_{a + 1} \rightarrow \mathbb{R}$ and $\nu > 0$. The $\nu^{\text{th}}$-order nabla sum of $u$ based at $a$ is given by
\begin{equation}
\nonumber \big{(}\nabla ^{-\nu}_{a}u\big{)}(t) = \sum^{t}_{s = a + 1}H_{\nu - 1}(t, \rho(s))u(s), \quad t \in \mathbb{N}_{a},
\end{equation}
where by convention $\big{(}\nabla ^{-\nu}_{a}u\big{)}(a) = 0$.
\end{definition}

\begin{definition} \cite{Go}
Let $u: \mathbb{N}_{a + 1} \rightarrow \mathbb{R}$ and $0 < \nu \leq 1$. The $\nu^{\text{th}}$-order {\it Riemann--Liouville} nabla difference of $u$ based at $a$ is given by
\begin{equation}
\nonumber \big{(}\nabla ^{\nu}_{a}u\big{)}(t) = \Big{(}\nabla \big{(}\nabla_{a}^{-(1 - \nu)}u\big{)}\Big{)}(t), \quad t\in\mathbb{N}_{a + 1}.
\end{equation}
\end{definition}

Ahrendt et al. \cite{Ah} showed that the definition of a fractional difference can be rewritten in a form similar to the definition of a fractional sum.

\begin{theorem} \cite{Ah} \label{Ah}
Let $u: \mathbb{N}_{a + 1} \rightarrow \mathbb{R}$ and $0 < \nu < 1$. Then,
\begin{equation}
\nonumber \big{(}\nabla ^{\nu}_{a}u\big{)}(t) = \sum^{t}_{s = a + 1}H_{-\nu - 1}(t, \rho(s))u(s), \quad t\in\mathbb{N}_{a + 1}.
\end{equation}
\end{theorem}

\begin{definition} \cite{An}
Let $u: \mathbb{N}_{a} \rightarrow \mathbb{R}$ and $0 < \nu \leq 1$. The $\nu^{\text{th}}$-order {\it Caputo} nabla fractional difference of $u$ based at $a$ is given by
\begin{equation}
\nonumber \big{(}\nabla ^{\nu}_{a*}u\big{)}(t) = \Big{(}\nabla_{a}^{-(1 - \nu)} \big{(}\nabla u\big{)}\Big{)}(t), \quad t\in\mathbb{N}_{a + 1}.
\end{equation}
\end{definition}

The following identity is useful in transforming the Caputo nabla fractional difference into the Riemann--Liouville nabla fractional difference.

\begin{theorem} \cite{Ab 1} \label{Relation} 
Let $u: \mathbb{N}_{a} \rightarrow \mathbb{R}$ and $0 < \nu < 1$. Then, 
\begin{equation}
\nonumber \big{(}\nabla ^{\nu}_{a*}u\big{)}(t) = \big{(}\nabla ^{\nu}_{a}u\big{)}(t) - H_{-\nu}(t, a)u(a), \quad t\in\mathbb{N}_{a + 1}.
\end{equation}
\end{theorem}

\section{Exponential Functions of Nabla Fractional Calculus}
Acar et al. \cite{Ac} and Nagai \cite{Na} introduced the exponential functions of nabla fractional calculus as the unique solutions of the following initial value problems associated with the Riemann--Liouville and the Caputo nabla fractional differences: 
\begin{equation} \label{FDE RL}
\begin{cases}
\big{(}\nabla^{\nu}_{\rho(0)}w\big{)}(t) = \lambda w(t), \quad t \in \mathbb{N}_{1}, \\
\big{(}\nabla^{- (1 - \nu)}_{\rho(0)}w\big{)}(0) = w(0) = 1,
\end{cases}
\end{equation}
and 
\begin{equation} \label{FDE C}
\begin{cases}
\big{(}\nabla^{\nu}_{0*}x\big{)}(t) = \lambda x(t), \quad t \in \mathbb{N}_{1}, \\
x(0) = 1,
\end{cases}
\end{equation}
where $0 < \nu < 1$ and $\left|\lambda\right| < 1$. The unique solutions of the initial value problems \eqref{FDE RL} and \eqref{FDE C} are represented by $\hat{e}_{\nu, \nu}(\lambda, t^{\overline{\nu}})$ and $\hat{e}_{\nu}(\lambda, t^{\overline{\nu}})$, respectively, where 
\begin{equation} \label{E RL}
\hat{e}_{\nu, \nu}(\lambda, t^{\overline{\nu}}) = \sum^{\infty}_{k = 0} \lambda^k H_{\nu k + \nu - 1}(t, \rho(0)), \quad t \in \mathbb{N}_{0},
\end{equation}
and 
\begin{equation} \label{E C}
\hat{e}_{\nu}(\lambda, t^{\overline{\nu}}) = \sum^{\infty}_{k = 0} \lambda^k H_{\nu k}(t, 0), \quad t \in \mathbb{N}_{0}.
\end{equation}
Atici et al. \cite{At 2}, \v{C}erm\'{a}k et al. \cite{Ce}, Eloe et al. \cite{Eloe}, Jia et al. \cite{Ji} and Wu et al. \cite{Wu} obtained the following asymptotic results of the discrete exponential functions. 
\begin{align}
\lim_{t \rightarrow \infty} \hat{e}_{\nu, \nu}(\lambda, t^{\overline{\nu}}) & = 0, \quad \lambda \in (-1, 0], \label{FDE RL A} \\
\lim_{t \rightarrow \infty} \hat{e}_{\nu, \nu}(\lambda, t^{\overline{\nu}}) & = \infty, \quad \lambda \in (0, 1), \label{FDE RL AA} \\
\lim_{t \rightarrow \infty} \hat{e}_{\nu}(\lambda, t^{\overline{\nu}}) & = 0, \quad \lambda \in (-1, 0), \label{FDE C A} \\
\lim_{t \rightarrow \infty} \hat{e}_{\nu}(\lambda, t^{\overline{\nu}}) & = \infty, \quad \lambda \in (0, 1). \label{FDE C AA}
\end{align}

Using triangular strip matrices, Podlubny \cite{Po 1} described a matrix approach to find numerical solutions of fractional differential equations. Motivated by this technique, we present a matrix method to compute the exponential functions \eqref{E RL} and \eqref{E C} numerically.

\subsection{Computation of \eqref{E RL}:} Let $m \in \mathbb{N}_{1}$ and consider the initial value problem associated with \eqref{FDE RL}:
\begin{equation} \label{FDE RL 1}
\begin{cases}
\big{(}\nabla^{\nu}_{\rho(0)}w\big{)}(t) = \lambda w(t), \quad t \in \mathbb{N}^{m}_{1}, \\
\big{(}\nabla^{- (1 - \nu)}_{\rho(0)}w\big{)}(0) = w(0) = 1.
\end{cases}
\end{equation}
Rewriting the equation in \eqref{FDE RL 1} using Theorem \ref{Ah}, we have
\begin{equation} \label{FDE RL 11}
\sum^{t}_{s = 0}H_{-\nu - 1}(t, \rho(s))w(s) = \lambda w(t), \quad t \in \mathbb{N}^{m}_{1}.
\end{equation}
Rearranging the terms in \eqref{FDE RL 11}, we obtain
\begin{equation} \label{FDE RL 111}
(1 - \lambda) w(t) + \sum^{t - 1}_{s = 1}H_{-\nu - 1}(t, \rho(s))w(s) = - H_{-\nu - 1}(t, \rho(0))w(0), \quad t \in \mathbb{N}^{m}_{1}.
\end{equation}
Denote by $\tilde{w} = [w(1), w(2), \cdots, w(m)]^T$. Then, the matrix form of \eqref{FDE RL 111} is given by
\begin{equation}
\nonumber \mathcal{L}\tilde{w} = - \mathcal{B},
\end{equation}
where
\begin{equation}
\resizebox{.95 \hsize}{!}{$\nonumber \mathcal{L} = \left(
\begin{array}{ccccccc}
1 - \lambda & 0  & \cdots & \cdots & 0 & 0 \\
H_{-\nu - 1}(2, \rho(1)) & 1 - \lambda  & \cdots & \cdots & 0 & 0 \\
H_{-\nu - 1}(3, \rho(1)) & H_{-\nu - 1}(3, \rho(2))  & \cdots & \cdots & 0 & 0 \\
\vdots & \vdots & \vdots  & \vdots & \vdots & \vdots \\
\vdots & \vdots & \vdots  & \vdots & \vdots & \vdots \\
H_{-\nu - 1}(m - 1, \rho(1)) & H_{-\nu - 1}(m - 1, \rho(2)) & \cdots & \cdots & 1 - \lambda & 0 \\
H_{-\nu - 1}(m, \rho(1)) & H_{-\nu - 1}(m, \rho(2)) & \cdots & \cdots & H_{-\nu - 1}(m, \rho(m - 1)) & 1 - \lambda\\
\end{array}
\right)_{m \times m}$}
\end{equation}
is a lower triangular strip matrix and
\begin{equation}
\nonumber \mathcal{B} = \left(
      \begin{array}{c}
        H_{-\nu - 1}(1, \rho(0)) \\
        H_{-\nu - 1}(2, \rho(0)) \\
        H_{-\nu - 1}(3, \rho(0)) \\
        \vdots \\
        \vdots \\
        H_{-\nu - 1}(m - 1, \rho(0)) \\
        H_{-\nu - 1}(m, \rho(0)) \\
        \end{array}
    \right)_{m \times 1}.
\end{equation}
Since $\mathcal{L}$ is non-singular, the exponential function \eqref{E RL} can be computed by the following numerical algorithm: $$\hat{e}_{\nu, \nu}(\lambda, t^{\overline{\nu}}) = - \mathcal{L}^{-1}\mathcal{B}, \quad t \in \mathbb{N}^{m}_{1}.$$ Here $\mathcal{L} = \left[\mathcal{L}_{ij}\right]_{m \times m}$ and $\mathcal{B} = \left[\mathcal{B}_{i}\right]_{m \times 1}$, where
\begin{equation*}
\mathcal{L}_{ij} = 
\begin{cases}
1 - \lambda, \hspace{0.77 in} i = j, \\
0, \hspace{1.07 in} i < j, \\
H_{-\nu - 1}(i, \rho(j)), \hspace{0.15 in} i > j,
\end{cases}
\end{equation*}
and 
\begin{equation*}
\mathcal{B}_{i} = H_{-\nu - 1}(i, \rho(0)).
\end{equation*}

\begin{exe}
Computation of $\hat{e}_{0.5, 0.5}(-0.5, t^{\overline{0.5}})$ for $t \in \mathbb{N}^{10}_{1}$:
\end{exe}

We have
\begin{equation*}
\resizebox{.95 \hsize}{!}{$\nonumber \mathcal{L} = \left(\begin{array}{cccccccccc} 1.5000 & 0 & 0 & 0 & 0 & 0 & 0 & 0 & 0 & 0\\ -0.5000 & 1.5000 & 0 & 0 & 0 & 0 & 0 & 0 & 0 & 0\\ -0.1250 & -0.5000 & 1.5000 & 0 & 0 & 0 & 0 & 0 & 0 & 0\\ -0.0625 & -0.1250 & -0.5000 & 1.5000 & 0 & 0 & 0 & 0 & 0 & 0\\ -0.0391 & -0.0625 & -0.1250 & -0.5000 & 1.5000 & 0 & 0 & 0 & 0 & 0\\ -0.0273 & -0.0391 & -0.0625 & -0.1250 & -0.5000 & 1.5000 & 0 & 0 & 0 & 0\\ -0.0205 & -0.0273 & -0.0391 & -0.0625 & -0.1250 & -0.5000 & 1.5000 & 0 & 0 & 0\\ -0.0161 & -0.0205 & -0.0273 & -0.0391 & -0.0625 & -0.1250 & -0.5000 & 1.5000 & 0 & 0\\ -0.0131 & -0.0161 & -0.0205 & -0.0273 & -0.0391 & -0.0625 & -0.1250 & -0.5000 & 1.5000 & 0\\ -0.0109 & -0.0131 & -0.0161 & -0.0205 & -0.0273 & -0.0391 & -0.0625 & -0.1250 & -0.5000 & 1.5000 \end{array}\right)$},
\end{equation*}
\begin{equation*}
\nonumber \mathcal{B} = \left(\begin{array}{c} -0.5000\\ -0.1250\\ -0.0625\\ -0.0391\\ -0.0273\\ -0.0205\\ -0.0161\\ -0.0131\\ -0.0109\\ -0.0093 \end{array}\right).
\end{equation*}
Then, for $t \in \mathbb{N}^{10}_{1}$,
\begin{equation*}
\nonumber \hat{e}_{0.5, 0.5}(-0.5, t^{\overline{0.5}}) = - \mathcal{L}^{-1}\mathcal{B} = \left(\begin{array}{c} 0.3333\\ 0.1944\\ 0.1343\\ 0.1009\\ 0.0798\\ 0.0654\\ 0.0550\\ 0.0472\\ 0.0411\\ 0.0362 \end{array}\right).
\end{equation*}

\begin{exe}
Computation of $\hat{e}_{0.5, 0.5}(0.5, t^{\overline{0.5}})$ for $t \in \mathbb{N}^{10}_{1}$:
\end{exe}

We have
\begin{equation*}
\resizebox{.95 \hsize}{!}{$\nonumber \mathcal{L} = \left(\begin{array}{cccccccccc} 0.5000 & 0 & 0 & 0 & 0 & 0 & 0 & 0 & 0 & 0\\ -0.5000 & 0.5000 & 0 & 0 & 0 & 0 & 0 & 0 & 0 & 0\\ -0.1250 & -0.5000 & 0.5000 & 0 & 0 & 0 & 0 & 0 & 0 & 0\\ -0.0625 & -0.1250 & -0.5000 & 0.5000 & 0 & 0 & 0 & 0 & 0 & 0\\ -0.0391 & -0.0625 & -0.1250 & -0.5000 & 0.5000 & 0 & 0 & 0 & 0 & 0\\ -0.0273 & -0.0391 & -0.0625 & -0.1250 & -0.5000 & 0.5000 & 0 & 0 & 0 & 0\\ -0.0205 & -0.0273 & -0.0391 & -0.0625 & -0.1250 & -0.5000 & 0.5000 & 0 & 0 & 0\\ -0.0161 & -0.0205 & -0.0273 & -0.0391 & -0.0625 & -0.1250 & -0.5000 & 0.5000 & 0 & 0\\ -0.0131 & -0.0161 & -0.0205 & -0.0273 & -0.0391 & -0.0625 & -0.1250 & -0.5000 & 0.5000 & 0\\ -0.0109 & -0.0131 & -0.0161 & -0.0205 & -0.0273 & -0.0391 & -0.0625 & -0.1250 & -0.5000 & 0.5000 \end{array}\right)$},
\end{equation*}
\begin{equation*}
\nonumber \mathcal{B} = \left(\begin{array}{c} -0.5000\\ -0.1250\\ -0.0625\\ -0.0391\\ -0.0273\\ -0.0205\\ -0.0161\\ -0.0131\\ -0.0109\\ -0.0093 \end{array}\right).
\end{equation*}
Then, for $t \in \mathbb{N}^{10}_{1}$,
\begin{equation*}
\nonumber \hat{e}_{0.5, 0.5}(0.5, t^{\overline{0.5}}) = - \mathcal{L}^{-1}\mathcal{B} = \left(\begin{array}{c} 1\\ 1.2500\\ 1.6250\\ 2.1406\\ 2.8359\\ 3.7676\\ 5.0127\\ 6.6749\\ 8.8925\\ 11.8505 \end{array}\right).
\end{equation*}

\begin{exe}
The graphs of $\hat{e}_{0.5, 0.5}(-0.5, t^{\overline{0.5}})$ and $\hat{e}_{0.5, 0.5}(0.5, t^{\overline{0.5}})$ for $t \in \mathbb{N}^{100}_{1}$ are shown in Figures 1 and 2, respectively.
\end{exe}

\begin{center}
\begin{figure}[!ht]
\includegraphics[width = 5.00 in]{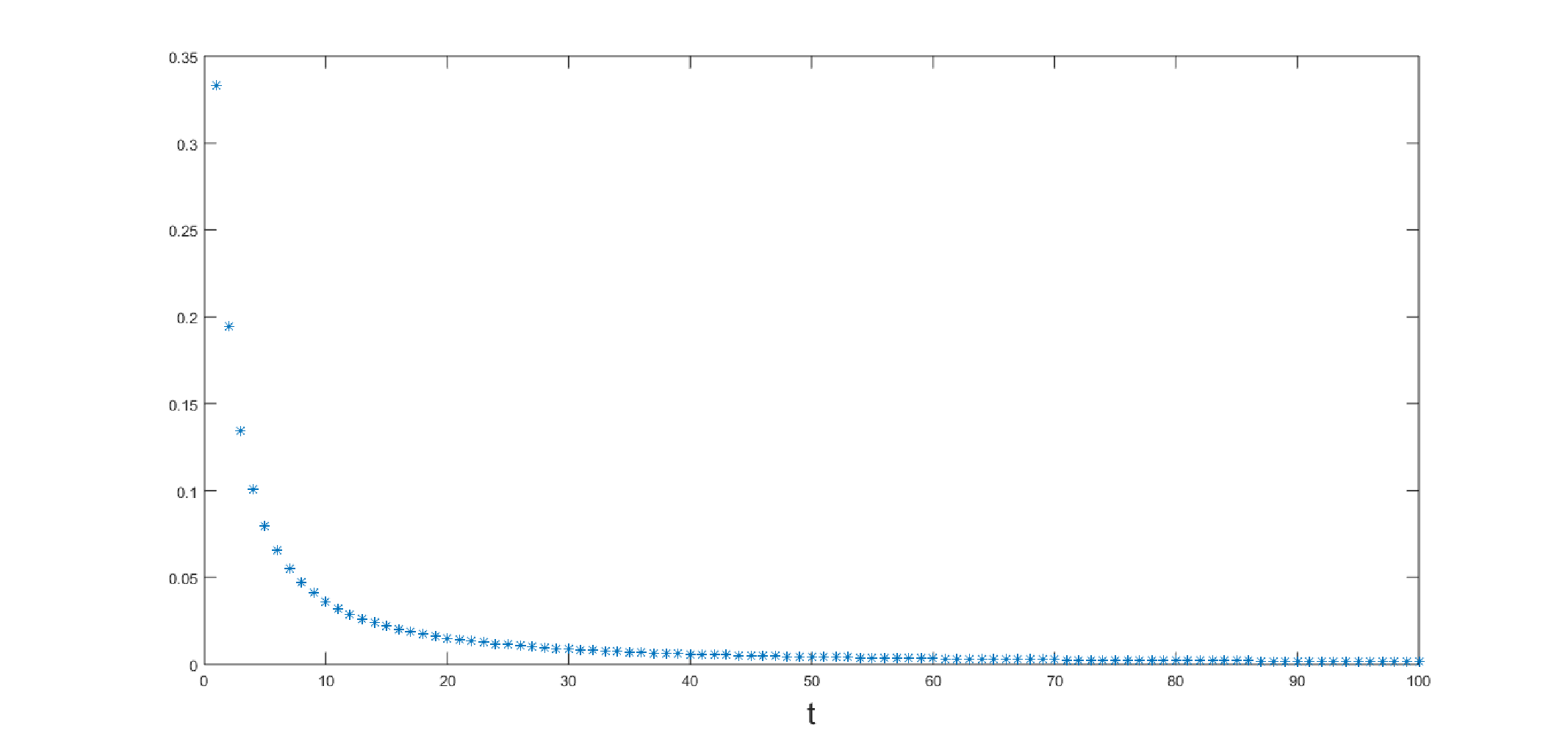}\\
\caption{}
\label{Figure 1}
\end{figure}
\end{center}

\begin{center}
\begin{figure}[!ht]
\includegraphics[width = 5.00 in]{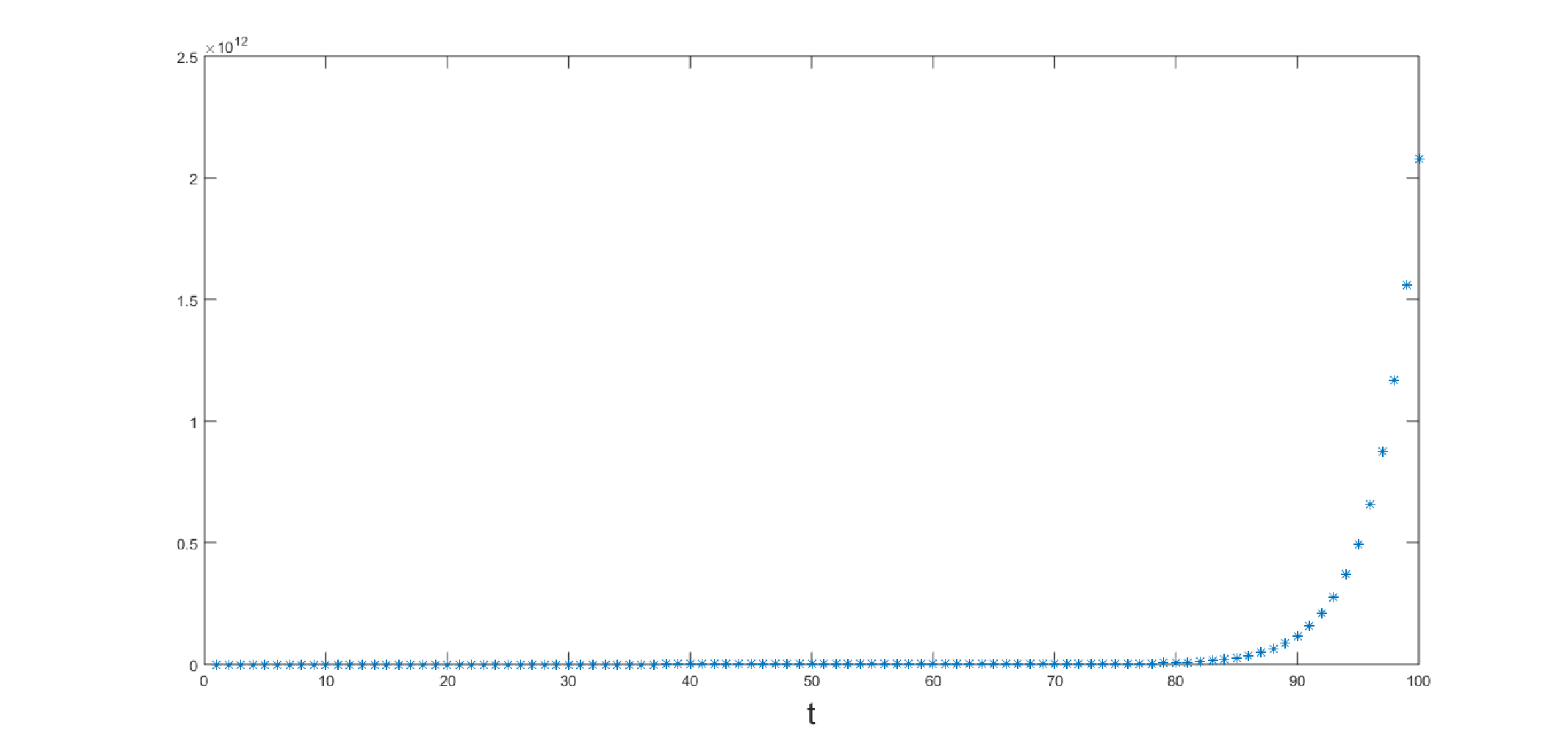}\\
\caption{}
\label{Figure 2}
\end{figure}
\end{center}

\subsection{Computation of \eqref{E C}:} Let $m \in \mathbb{N}_{1}$ and consider the initial value problem associated with \eqref{FDE C}:
\begin{equation} \label{FDE C 1}
\begin{cases}
\big{(}\nabla^{\nu}_{0*}x\big{)}(t) = \lambda x(t), \quad t \in \mathbb{N}^m_{1}, \\
x(0) = 1.
\end{cases}
\end{equation}
Rewriting the equation in \eqref{FDE C 1} using Theorem \ref{Ah} and Theorem \ref{Relation}, we have
\begin{equation} \label{FDE C 11}
\sum^{t}_{s = 1}H_{-\nu - 1}(t, \rho(s))x(s) - H_{-\nu}(t, 0)x(0) = \lambda x(t), \quad t \in \mathbb{N}^{m}_{1}.
\end{equation}
Rearranging the terms in \eqref{FDE C 11}, we obtain
\begin{equation} \label{FDE C 111}
(1 - \lambda) x(t) + \sum^{t - 1}_{s = 1}H_{-\nu - 1}(t, \rho(s))x(s) = H_{-\nu}(t, 0)x(0), \quad t \in \mathbb{N}^{m}_{1}.
\end{equation}
Denote by $\tilde{x} = [x(1), x(2), \cdots, x(m)]^T$. Then, the matrix form of \eqref{FDE C 111} is given by
\begin{equation}
\nonumber \mathcal{L}\tilde{x} = \mathcal{C},
\end{equation}
where
\begin{equation}
\nonumber \mathcal{C} = \left(
      \begin{array}{c}
        H_{-\nu}(1, 0) \\
        H_{-\nu}(2, 0) \\
        H_{-\nu}(3, 0) \\
        \vdots \\
        \vdots \\
        H_{-\nu}(m - 1, 0) \\
        H_{-\nu}(m, 0) \\
        \end{array}
    \right)_{m \times 1}.
\end{equation}
Since $\mathcal{L}$ is non-singular, the exponential function \eqref{E C} can be computed by the following numerical algorithm: $$\hat{e}_{\nu}(\lambda, t^{\overline{\nu}}) = \mathcal{L}^{-1}\mathcal{C}, \quad t \in \mathbb{N}^{m}_{1}.$$ Here $\mathcal{L} = \left[\mathcal{L}_{ij}\right]_{m \times m}$ and $\mathcal{C} = \left[\mathcal{C}_{i}\right]_{m \times 1}$, where
\begin{equation*}
\mathcal{L}_{ij} = 
\begin{cases}
1 - \lambda, \hspace{0.77 in} i = j, \\
0, \hspace{1.07 in} i < j, \\
H_{-\nu - 1}(i, \rho(j)), \hspace{0.15 in} i > j,
\end{cases}
\end{equation*}
and 
\begin{equation*}
\mathcal{C}_{i} = H_{-\nu}(i, 0).
\end{equation*}

\begin{exe}
Computation of $\hat{e}_{0.5}(-0.5, t^{\overline{0.5}})$ for $t \in \mathbb{N}^{10}_{1}$:
\end{exe}

We have
\begin{equation*}
\nonumber \mathcal{C} = \left(\begin{array}{c} 1\\ 0.5000\\ 0.3750\\ 0.3125\\ 0.2734\\ 0.2461\\ 0.2256\\ 0.2095\\ 0.1964\\ 0.1855 \end{array}\right).
\end{equation*}
Then, from Example 1, for $t \in \mathbb{N}^{10}_{1}$,
\begin{equation*}
\nonumber \hat{e}_{0.5}(-0.5, t^{\overline{0.5}}) = \mathcal{L}^{-1}\mathcal{C} = \left(\begin{array}{c} 0.6667\\ 0.5556\\ 0.4907\\ 0.4460\\ 0.4124\\ 0.3857\\ 0.3639\\ 0.3456\\ 0.3299\\ 0.3162 \end{array}\right).
\end{equation*}

\begin{exe}
Computation of $\hat{e}_{0.5}(0.5, t^{\overline{0.5}})$ for $t \in \mathbb{N}^{10}_{1}$:
\end{exe}

We have
\begin{equation*}
\nonumber \mathcal{C} = \left(\begin{array}{c} 1\\ 0.5000\\ 0.3750\\ 0.3125\\ 0.2734\\ 0.2461\\ 0.2256\\ 0.2095\\ 0.1964\\ 0.1855 \end{array}\right).
\end{equation*}
Then, from Example 2, for $t \in \mathbb{N}^{10}_{1}$,
\begin{equation*}
\nonumber \hat{e}_{0.5}(0.5, t^{\overline{0.5}}) = \mathcal{L}^{-1}\mathcal{C} = \left(\begin{array}{c} 2\\ 3\\ 4.2500\\ 5.8750\\ 8.0156\\ 10.8516\\ 14.6191\\ 19.6318\\ 26.3067\\ 35.1992 \end{array}\right).
\end{equation*}

\begin{exe}
The graphs of $\hat{e}_{0.5}(0.5, t^{\overline{0.5}})$ and $\hat{e}_{0.5}(-0.5, t^{\overline{0.5}})$ for $t \in \mathbb{N}^{100}_{1}$ are shown in Figures 3 and 4, respectively.
\end{exe}

\begin{center}
\begin{figure}[!ht]
\includegraphics[width = 5.00 in]{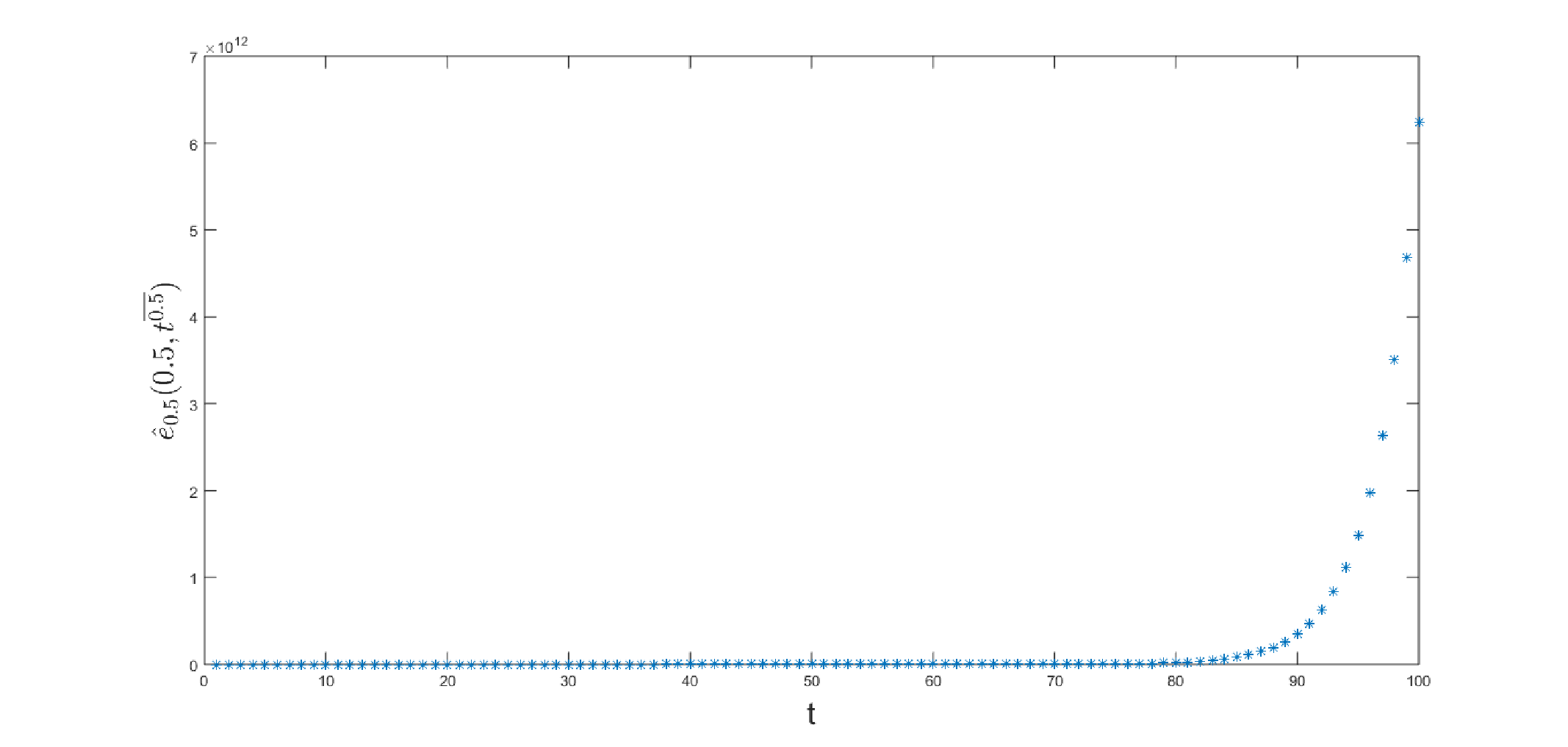}\\
\caption{}
\label{Figure 3}
\end{figure}
\end{center}

\begin{center}
\begin{figure}[!ht]
\includegraphics[width = 5.00 in]{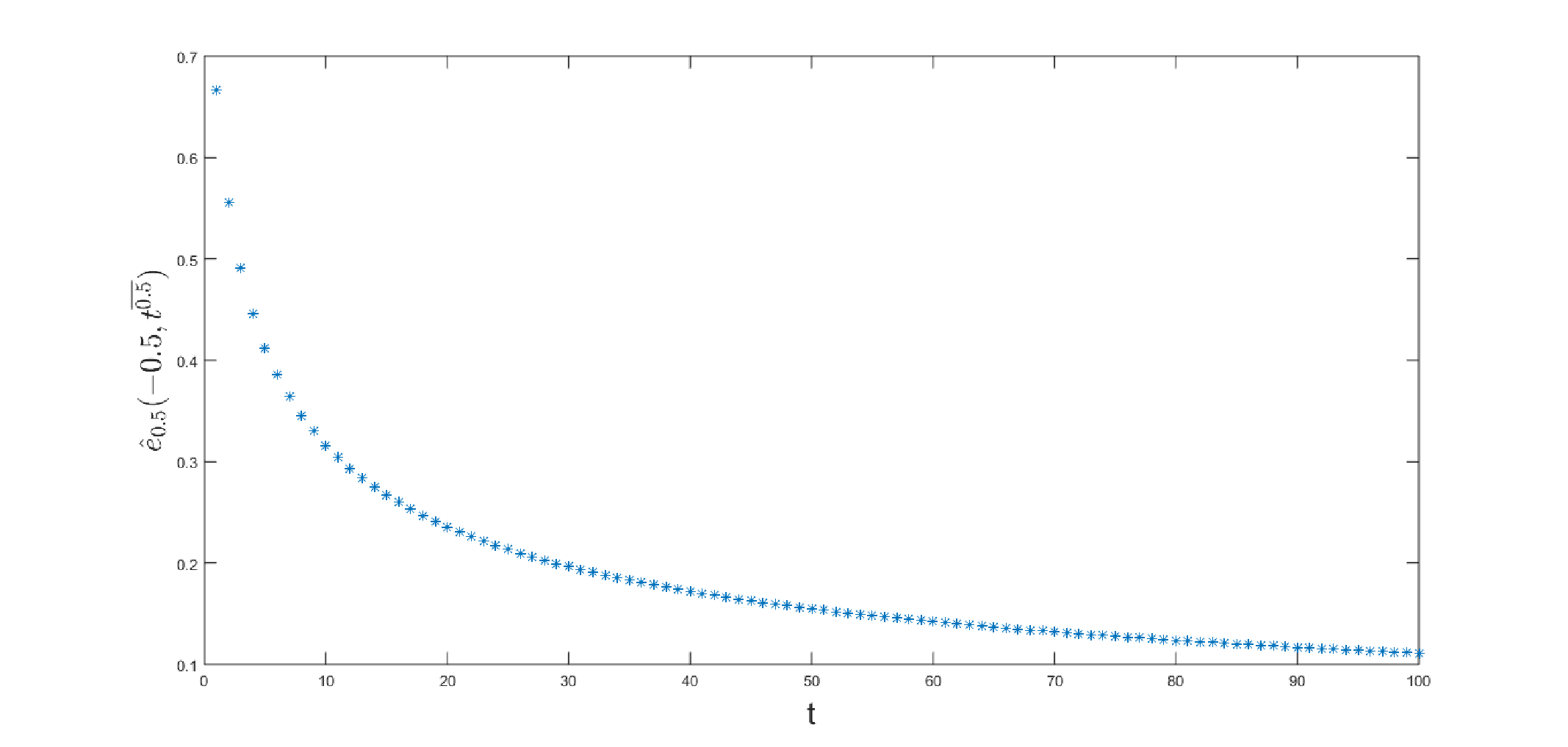}\\
\caption{}
\label{Figure 4}
\end{figure}
\end{center}

\section{Extensions} The method described in Section 2 can be extended to obtain numerical solutions of initial value problems involving linear non-homogeneous nabla fractional difference equations.

Let $0 < \nu < 1$ and $m \in \mathbb{N}_{1}$. Consider the initial value problem
\begin{equation} \label{IVP 1}
\begin{cases}
\big{(}\nabla^{\nu}_{\rho(0)} u\big{)}(t) = a(t) u(t) + f(t), \quad t \in \mathbb{N}^{m}_{1}, \\
\big{(}\nabla^{- (1 - \nu)}_{\rho(0)} u\big{)}(0) = u(0) = c,
\end{cases}
\end{equation}
where $a$, $f : \mathbb{N}^{m}_{1} \rightarrow \mathbb{R}$ such that $$a(t) \neq 1, \quad t \in \mathbb{N}^{m}_{1}.$$ Denote by $\tilde{u} = [u(1), u(2), \cdots, u(m)]^T$ and $\mathcal{F} = [f(1), f(2), \cdots, f(m)]^T$. Then, the matrix form of \eqref{IVP 1} is given by
\begin{equation}
\nonumber \mathcal{L}_1\tilde{u} = \mathcal{F} - c\mathcal{B},
\end{equation}
where
\begin{equation}
\resizebox{.95 \hsize}{!}{$\nonumber \mathcal{L}_1 = \left(
\begin{array}{ccccccc}
1 - a(1) & 0  & \cdots & \cdots & 0 & 0 \\
H_{-\nu - 1}(2, \rho(1)) & 1 - a(2)  & \cdots & \cdots & 0 & 0 \\
H_{-\nu - 1}(3, \rho(1)) & H_{-\nu - 1}(3, \rho(2))  & \cdots & \cdots & 0 & 0 \\
\vdots & \vdots & \vdots  & \vdots & \vdots & \vdots \\
\vdots & \vdots & \vdots  & \vdots & \vdots & \vdots \\
H_{-\nu - 1}(m - 1, \rho(1)) & H_{-\nu - 1}(m - 1, \rho(2)) & \cdots & \cdots & 1 - a(m - 1) & 0 \\
H_{-\nu - 1}(m, \rho(1)) & H_{-\nu - 1}(m, \rho(2)) & \cdots & \cdots & H_{-\nu - 1}(m, \rho(m - 1)) & 1 - a(m)\\
\end{array}
\right)_{m \times m}$}
\end{equation}
is a lower triangular strip matrix. Since $\mathcal{L}_1$ is non-singular, the solution of \eqref{IVP 1} can be computed by the following numerical algorithm: $$u(t) = \mathcal{L}_1^{-1}[\mathcal{F} - c\mathcal{B}], \quad t \in \mathbb{N}^{m}_{1}.$$ 

Replacing the $\nu$-th order Riemann--Liouville nabla fractional difference operator $\nabla^{\nu}_{\rho(0)}$ in \eqref{IVP 1} with the $\nu$-th order Caputo operator $\nabla^{\nu}_{0*}$, the matrix form of the initial value problem
\begin{equation} \label{IVP 2}
\begin{cases}
\big{(}\nabla^{\nu}_{0*} u\big{)}(t) = a(t) u(t) + f(t), \quad t \in \mathbb{N}^{m}_{1}, \\
u(0) = c,
\end{cases}
\end{equation}
is given by
\begin{equation}
\nonumber \mathcal{L}_1\tilde{u} = \mathcal{F} + c\mathcal{C}.
\end{equation}
Since $\mathcal{L}_1$ is non-singular, the numerical solution of \eqref{IVP 2} can be computed by the following numerical algorithm: $$u(t) = \mathcal{L}_1^{-1}[\mathcal{F} + c\mathcal{C}], \quad t \in \mathbb{N}^{m}_{1}.$$

\end{document}